\documentclass[a4paper,11pt]{article}
\usepackage{graphicx,amsmath,amsfonts,amssymb,amsthm,psfrag,color}
\usepackage{bbm,wrapfig,float,framed}
\usepackage{fullpage}

\newcommand{\ds}{\displaystyle}
\newcommand{\dint}{\mathrm{d}}

\newcommand{\eps}{\varepsilon}

\newtheorem{thm}{Theorem}[section]

\numberwithin{equation}{section}

\begin{document}
\title{The walk on moving spheres: a new tool
for simulating Brownian motion's exit time from a domain}
\author{M. Deaconu\\
Inria, Villers-l{\`e}s-Nancy, F-54600, France\\
Universit{\'e} de Lorraine, CNRS, Institut Elie Cartan de Lorraine - UMR 7502,\\
 Vandoeuvre-l{\`e}s-Nancy, F-54506, France\\
 \texttt{Madalina.Deaconu@inria.fr} \and S. Herrmann \\
 IMB UMR5584, CNRS, Univ. Bourgogne Franche-Comt\'e,\\ F-21000 Dijon, France\\
 \texttt{Samuel.Herrmann@u-bourgogne.fr} \and S. Maire\\
 Aix-Marseille Universit\'e, CNRS, LSIS, UMR 7296, F-13397 Marseille France \\
Universit\'e de Toulon CNRS, LSIS, UMR 7296, F-83957 La Garde France\\
\texttt{maire@univ-tln.fr}}
\date{\today}
\maketitle
\begin{abstract}
In this paper we introduce a new method for the simulation of the exit time and  exit position of a $\delta$-dimensional Brownian motion from a domain. The main interest of our method is that it avoids splitting time schemes as well as inversion of complicated series.  The method, called walk on moving spheres algorithm, was first introduced for hitting times of Bessel processes. In this study this method is adapted and developed for the first time for the Brownian motion hitting times.  The idea is to use the connexion between the $\delta$-dimensional Bessel process and the $\delta$-dimensional Brownian motion thanks to an explicit Bessel hitting time distribution associated with a particular curved boundary. This allows to build a fast and accurate numerical scheme for approximating the hitting time. We introduce also an overview of existing methods for the simulation of the Brownian hitting time and perform numerical comparisons with existing methods.  
\end{abstract}
%\begin{keyword} Walk on moving spheres method \sep 
%Bessel processes \sep Brownian hitting time \sep numerical algorithm
%\end{keyword}
%\end{frontmatter}

%%%%%%%%%%%%%%%%%%%%%%%%%%%%%%%%%%
\section{Introduction}
%%%%%%%%%%%%%%%%%%%%%%%%%%%%%%%%%%
Computing the first hitting time of a boundary by a stochastic process with a high accuracy is of great interest for many areas of applications. Examples range from neuronal sciences, financial derivatives with barriers, optimal stopping problems and so on. For general stochastic diffusion processes, the simulation of the exit time from a domain is in general obtained by the Euler scheme. While naive versions of this scheme reach an order one half for the computations of weak approximations, it is possible to obtain order one approximations thanks to a barrier correction \cite{gobet01}. The number of steps before hitting the boundary is nevertheless  proportional to the inverse of the time discretisation step. However, when the diffusion process reduces to a standard multidimensional Brownian motion, alternative more efficient simulation methods can be used.   The random walk on spheres (WOS) introduced by Muller \cite{muller56} relies on the isotropy of the Brownian motion and enables to make large jumps instead of small ones for the Euler scheme. Its mean number of steps before hitting the boundary is proportional to  $\left|\ln(\varepsilon)\right|$ where $\varepsilon$ is the parameter of the absorption boundary layer.   The  random walk on rectangles  proposed by Deaconu and Lejay \cite{deaconu-lejay08} uses the same ideas and may be even more efficient for a polygonal domain. Using these two methods, the elapsed time is nevertheless a lot harder to simulate than with the Euler scheme. Indeed its simulation needs the inversion method on a cumulative distribution function which is a complicated series. A fast and accurate simulation for the hitting time law is mandatory for the computation of, for instance, the principal eigenvalue of the Laplace operator \cite{lejay-maire07}.  We describe and study here a method called the walk on  moving spheres (WOMS), introduced by Deaconu and Herrmann \cite{deaconu-herrmann13} which conciliates a small number of steps before absorption and an easy way to simulate the exit time. 

The paper is organised as follows. In section two, we recall the properties of the WOS method and discuss the exit time of a sphere. The study of the hitting time and exit position methods for a Brownian motion is performed in section 3 and the Bessel hitting time is also introduced. The section 4 is devoted to the new method based on the simulation of the Brownian position by an uniform random variable and the hitting time by using the explicit expression of the Bessel hitting time distribution.  The last section illustrates numerical results and compares the mean number of steps before absorption and the efficiency of the WOS and WOMS methods for exit time simulation on a simple  numerical example. 
%\\[5pt]
%%%%%%%%%%%%%%%%%%%%%%%%%%%%%%%%
\section{Random walk on spheres}
%%%%%%%%%%%%%%%%%%%%%%%%%%%%%%%%
\label{sec:muller}
The study  of the hitting time of a given boundary for the Brownian motion is of great interest in many applications. This research has a long history as it is connected  with the solution $u$ of the Dirichlet problem :
\begin{equation}
\label{dirichlet}
\left\{
\begin{array}{ll}
\ds\frac{1}{2}\Delta u(x)& = 0 \mbox { on } D,\\
u(x)& = f(x) \mbox { on } \partial D,
\end{array}
\right.
\end{equation}
where $D$ denotes a bounded finitely connected domain in $\mathbb{R}^\delta$ and $\partial D$ its boundary, assumed throughout this paper to be of sufficient regularity, in order to ensure that the Dirichlet problem has a unique solution. The function  $f$ is continuous on the boundary $\partial D$. The probabilistic approach of this problem is a powerful tool that allows to express the solution of $(\ref{dirichlet})$ in the form $u(x) = \mathbb{E}_x [f(B^{\delta,x}_{\tau_{\partial D}})]$ where $B_t^{\delta,x}$ stands for the $\delta$-dimensional Brownian motion starting from $x$,  and $\tau_{\partial D}$ its first hitting time of the boundary $\partial D$, that is
$\tau_{\partial D}=\inf\{ t>0: B_t^{\delta,x} \in \partial D\}.$

Introducing an efficient numerical approximation of the quantities $\tau_{\partial D}$ and $B_{\tau_{\partial D}}^{\delta,x}$ also provides an accurate procedure in order to approximate the solution to the Dirichlet problem (\ref{dirichlet}).

For these purposes, the methods that can be considered are mainly based on splitting time methods like the Euler scheme but they are not computationally efficient and may overestimate the hitting time. An alternative approach can use the explicit form of the distribution of the hitting time and uses the inversion method. However, this procedure involves  complicated series and special functions like the Bessel ones. 

One of the revolutionary ideas on this topic is due to Muller \cite{muller56} who introduced a method called the random walk on spheres (WOS). His method approaches the hitting time and the exit position for the Brownian motion starting from $x$ and living in the domain $D$. This procedure is based on Monte Carlo methods for solving Dirichlet problem. 
The idea of the algorithm is to start by constructing the largest sphere centred at $x$ and included in $D$. For the Brownian motion starting at $x$, we consider the hitting position of this sphere by choosing uniformly a point on its  boundary. This gives the new starting point and the new center of the largest sphere included in $D$, used for the second step of the algorithm. The algorithm generates then iteratively the first exit time and position for the Brownian motion, starting in the current point, from the largest sphere included in $D$ and centred at the current point. The algorithm stops when the exit position is as close as suited to the boundary $\partial D$. This method relies on analytical expressions of the distribution functions for the first exit time and the first exit position from a sphere which is the uniform law on this sphere.

Since then, the WOS method has been extended for many applications as for example in Sabelfeld and Talay \cite{sabelfeld-talay95}, 
%Hwang, Mascagni and Given, \cite{hwang-mascagni-given03} 
and Golyandina \cite{golyandina04}. Further generalisations of the method for non-homogeneous media are introduced by Milstein and coauthors \cite{milstein-rybkina93, milstein-tretyakov99}. For polygonal domains a similar method, based on random walk on rectangles, was proposed by Deaconu and Lejay \cite{deaconu-lejay08}.

In order to evaluate the exit position this procedure is really efficient. However, when considering the exit time, at each step of the algorithm one needs to evaluate the quantity 
\begin{equation}
\label{exit-time-mb}
\tau_L = \inf\{t>0 : \|B_t^{\delta,x}\| =L\}
\end{equation}
where $L$ is the radius of the corresponding sphere in the algorithm. This is the first time that the Euclidean norm of a $\delta$-dimensional Brownian motion hits the level $L$, and represents also the hitting time of the level $L$ for the $\delta$-dimensional Bessel process. Up to now, there is no general analytical formula allowing the numerical simulation of the distribution of $\tau_L$.  

For the case of the Bessel process starting from $x$ an explicit form of the Laplace transform of $\tau_L$ exists \cite{deaconu-herrmann13} for $x>0$:
\begin{equation}
\label{laplacebess}
 \mathbb{E}_{x}\Big[ e^{-\lambda \tau_L} \Big]=\frac{x^{-\nu}}{L^{-\nu}}\frac{I_\nu(x\sqrt{2\lambda})}{I_\nu(L\sqrt{2\lambda})}\ \mbox{and}\  \mathbb{E}_0\Big[ e^{-\lambda \tau_L} \Big]=\frac{(L\sqrt{2\lambda})^\nu}{2^\nu\Gamma(\nu+1)}\frac{1}{I_\nu(L\sqrt{2\lambda})}
\end{equation}
here $I_\nu(x)$ denotes the modified Bessel function and $\nu$ the index of the Bessel process, defined by $\nu=\delta/2-1$. Ciesielsky and Taylor \cite{Ciesielski-Taylor} proved that, for $\delta \in \mathbb{N}$ and $x=0$, the tail distribution is given by  
\begin{equation*}
 \mathbb{P}_0(\tau_L>t)=\frac{1}{2^{\nu-1}\Gamma(\nu+1)}\,\sum_{k=1}^{\infty}\frac{j_{\nu,k}^{\nu-1}}{\mathcal{J}_{\nu+1}(j_{\nu,k})}\,e^{-\frac{j_{\nu,k}^2}{2L^2}t},
%\quad x>0,
\end{equation*}
where $\mathcal{J}_\cdot$  is the Bessel function of the first kind, and $j_{\cdot,k}$ the associated sequence of its positive zeros. A similar formula is available for $x>0$.

Despite this explicit form, these formulas are obviously miss-adapted and not well suited for numerical purposes. We present in the next two sections the main results used for Brownian and Bessel cases. 
%%%%%%%%%%%%%%%%%%%%%%%%%%%%%%%%
\section{Hitting time of one-sided moving boundaries}
\label{sec:hitting}
\subsection{The one-dimensional Brownian case}
%%%%%%%%%%%%%%%%%%%%%%%%%%%%%%%
We shall first study the Brownian hitting time and afterwards focus on the Bessel case. The standard one-dimensional Brownian motion $(B_t,\ t\ge 0)$ satisfies nice hitting time properties and an explicit expression of the hitting time distribution for straight line boundaries is available.\\
Let $\psi:\mathbb{R}_+\to \mathbb{R}$ be a continuous function. We denote by $\tau_\psi=\inf\{t>0:\ B_t=\psi(t)\}$,
the first hitting time of the curved boundary $\psi$ for the Brownian motion.
\subsubsection{Hitting a given level $L$} 
Let $\psi$ be the constant function equal to $L>0$. Introduce the exponential martingale associated to the Brownian motion and use the optional stopping theorem to obtain, for any $\lambda\in\mathbb{R}$,
\[
\mathbb{E}[e^{-\lambda^2\tau_\psi/2}]=e^{-\lambda L}\mathbb{E}[e^{\lambda B_{\tau_\psi}-\lambda^2\tau_\psi/2}]=e^{-\lambda L}.
\]
In other words, 
\(
\mathbb{E}[e^{-\lambda\tau_\psi/2}]=e^{-L\sqrt{2\lambda}},\quad \lambda>0.
\)
In the exit time framework, we can often compute the Laplace transform of hitting times but inverting such expressions is not usually a simple task. For the present situation, the reflexion principle of the Brownian path and the scaling property permit to overcome this difficulty.  We have 
\begin{align*}
\mathbb{P}(\tau_\psi\le t)&=\mathbb{P}\Big(\sup_{0\le s \le t}B_s\ge L\Big)=\mathbb{P}(|B_t|\ge L)=\mathbb{P}\Big(|G|\ge \frac{L}{\sqrt{t}}\Big)=\mathbb{P}\Big(\frac{L^2}{G^2}\le t\Big),
\end{align*}
where $G$ stands for a standard normally distributed random variable. We deduce the identity
\begin{equation}
\label{eq:ident}
\tau_\psi\overset{\Delta}{=}\frac{L^2}{G^2}
\end{equation}
which immediately yields the probability density function (pdf) of $\tau_\psi$:
\begin{equation}
\label{eq:pdf}
p_\psi(t):=\frac{L}{\sqrt{2\pi}t^{3/2}}\, e^{-\frac{L^2}{2t}}.
\end{equation}
Let us introduce the function $\xi$ corresponding to the value of the Brownian pdf at the boundary, defined by:
\begin{equation}
\label{eq:defxi}
\xi(t)=\frac{1}{\sqrt{2\pi t}}\, \exp\left(-\frac{\psi^2(t)}{2t}\right),\quad t>0.
\end{equation}
An important feature is the relation between the pdf of the hitting time $p_\psi$ and the function $\xi$, as $\psi(t)=L$ for all $t$:
\begin{equation}
\label{eq:relat}
p_{\psi}(t)=b(t)\xi(t)\quad\mbox{with}\quad b(t)=\frac{L}{t}.
\end{equation}
To sum up, for the Brownian hitting time in the constant boundary case, the pdf has a simple expression and \eqref{eq:ident} is of prime interest for numerical purposes. 
\subsubsection{Hitting a straight line}
An explicit expression can also be obtained in the general straight line case. Let us assume now that 
\(
\psi(t)=L+\beta t, \quad \beta >0.
\)
\eqref{eq:relat} is still valid for a one sided boundary. The Bachelier-L{\'e}vy formula holds
\begin{equation}\label{eq:gene}
p_\psi(t)=\frac{L}{\sqrt{2\pi}t^{3/2}}\,\exp\left(-\frac{(L+\beta t)^2}{2t}\right),\quad t> 0.
\end{equation}
The proof relies on the Girsanov change of measure formula. For $\tilde{B}_t=B_t-\beta t$, we have $\tau_\psi(B)\overset{\Delta}{=}\tau_L(\tilde{B})$ where $\tau_L$ is the first hitting time of the level $L$ for $\tilde{B}$. Under the change of measure, $\tilde{B}_t$ becomes a Brownian motion. More precisely, defining $D_t=\exp\{-\beta B_t-\beta^2t/2 \}$, we have that 
%$\mathbb{P}(\tau_\psi(B)\le t)$ is equal to
\begin{align*}
\mathbb{P}(\tau_\psi(B)\le t)& 
=\mathbb{P}(\tau_L(\tilde{B})\le t)\\
& =\mathbb{E}\Big[1_{\{ \tau_L(B)\le t \}}D_{\tau_L}\Big] \\
&=\mathbb{E}\Big[ 1_{\{ \tau_L(B)\le t \}}e^{-{\frac{\beta^2\tau_L(B)}{2}}} \Big]e^{-L\beta}.
\end{align*}
We obtain thus \eqref{eq:gene} by using the time derivative and the explicit expression \eqref{eq:pdf} of the first passage time to the level $L$. The distribution arising here belongs to the inverse Gaussian family. More precisely $\tau_\psi$ has the inverse Gaussian distribution $I(-\frac{L}{\beta},L^2)$ (see for instance \cite{Devroye} p.148). 
%Consequently $(\beta\tau_\psi+a)^2/\tau_\psi$ is a chi-square distributed random variable with one degree of freedom. 
Consequently $\tau_\psi$ can be simulated with the simple generator introduced by Michael, Schucany and Haas \cite{Michael-76}. \\
Let us point out that in both preliminary cases (the constant boundary case and the straight line one), the particular relation between the hitting time pdf and the Brownian pdf at the boundary, given by \eqref{eq:relat}, is fulfilled with $b(t)=L/t$ and moreover the hitting times can be numerically easily generated.
\subsubsection{Hitting a general curved boundary: a numerical approach}
Obviously the general situation will not lead to simple pdf expressions like \eqref{eq:pdf} or \eqref{eq:gene}. Nevertheless Durbin \cite{Durbin-85,Durbin-92} proved that  \eqref{eq:relat} is a general formula with
\begin{equation}\label{eq:defdeb}
b(t)=\lim_{s\uparrow t}\frac{1}{t-s}\,\mathbb{E}\Big[  (\psi(s)-B_s)1_\Gamma\Big|B_t=\psi(t) \Big]\ \mbox{and}\ \Gamma:=\Big\{ \sup_{s\le u\le t}(B_u-\psi(u))\le 0 \Big\},
\end{equation}
as soon as the boundary is continuously differentiable. Even if the function $b$ is defined by a convergence result, it is often difficult to compute its value and its use usually requires an approximation procedure.\\ 
In \cite{Durbin-85}, the author computed the expression of $b(t)$ in the straight line case, using the formula \eqref{eq:defdeb}, and obtained as expected $b(t)=L/t$. Durbin noticed that the set $\Gamma$ does not play an important role in this particular situation. Indeed if the characteristic function is omitted in \eqref{eq:defdeb} the result is still valid:
\begin{equation}\label{eq:defdeb1}
b_1(t):=\lim_{s\uparrow t}\frac{1}{t-s}\,\mathbb{E}\left[ (\psi(s)-B_s)\Big|B_t=\psi(t)\right]=b(t)\quad \mbox{for}\quad \psi(t)=L+\beta t.
\end{equation}
This essential remark introduces a first rough approximation of the hitting time's pdf: $q_1(t)$ given by
\(
q_1(t):=b_1(t)\xi(t),
\)
where $b_1$ and $\xi$ are defined by \eqref{eq:defdeb1} and \eqref{eq:defxi}, respectively (this approximation is exact in the straight line case). It consists in fact in a tangent approximation of the boundary in a neighbourhood of $t$ suggested by Strassen \cite{strassen67}. By elementary Gaussian computations, we obtain
\(
b_1(t)
%\lim_{s\uparrow t}\frac{1}{t-s}\Big\{ \psi(s)-\mathbb{E}[B_s|B_t=\psi(t)] \Big\}=\lim_{s\uparrow t}\frac{1}{t-s}\Big\{ \psi(s)-\frac{\min(s,t)}{t}\,\psi(t) \Big\}
=\frac{\psi(t)}{t}-\psi'(t).
\)
Hence
\begin{equation}\label{eq:defdeq1}
q_1(t)=\left(\frac{\psi(t)}{t}-\psi'(t)\right) \xi(t)
\end{equation}
is the first approximation of $p_\psi(t)$. In order to get a sharper approximation, Durbin \cite{Durbin-92} proved that $p_\psi(t)$ solves a Volterra equation of the second type. In the Appendix of  \cite{Durbin-92}, Williams gave a more intuitive proof of this result. If $q(t,x,y)$ denotes the transition probabilities of the Brownian motion, then 
\begin{equation}
\label{eq:volterra}
p_\psi(t)=\left\{\frac{\psi(t)}{t}-\psi'(t)\right\}\,q(t,0,\psi(t))-\mathcal{P}_tp_\psi,
\end{equation}
where the operator $\mathcal{P}_t$ is defined by
\[
\mathcal{P}_tf=\int_0^tf(s)\left\{\frac{\psi(t)-\psi(s)}{t-s}-\psi'(t)\right\}\, q(t-s,\psi(s),\psi(t))\dint s.
\]
Observe that the first term in \eqref{eq:volterra} is exactly the approximation term $q_1(t)$. Consequently by developing a priori 
%$p_\psi(t)=q_1(t)-q_2(t)+\ldots+(-1)^{k-1}q_k(t)+\ldots$ leads to  
%\[
%q_1(t)-q_2(t)+\ldots+(-1)^{k-1}q_k(t)+\ldots=q_1(t)-\mathcal{P}_t\Big(q_1(t)-q_2(t)+\ldots+(-1)^{k-1}q_k(t)+\ldots\Big).
%\]
%We deduce that 
\(
p_\psi(t)=\sum_{k\ge 1}(-1)^{k-1}q_k(t)\)
we get $q_{k+1}(t)=\mathcal{P}_tq_k$ where $q_1$ is given by \eqref{eq:defdeq1}.
This discussion was made precise by Durbin who proposed an error bound for the approximation of the hitting time pdf by truncated series. This procedure permits to compute the pdf for any curved boundary using numerical integration. 
%in order to get the sequence $(q_k)_{k\ge 0}$. 
To sum up, the Brownian hitting time of a general curved boundary cannot be exactly described by an explicit expression of its pdf but can be generated via Durbin's approximation.
\subsubsection{Explicit expressions for particular curved boundaries}
This last paragraph concerning the Brownian hitting times, emphasises the use of the method of images developed by Lerche \cite{lerche_1986}. Since it is not possible to obtain nice expressions for a general boundary, we investigate families of boundaries which lead to explicit expressions of $p_\psi$.
% Of course in our study we take away the straight lines... 
The method of images is based on a positive, $\sigma$-finite measure $F$ (satisfying an integrability assumption) and a parameter $a>0$. The following function is then defined:
\[
h(t,x)=q(t,0,x)-\frac{1}{a}\int_0^\infty q(t,y,x)F(\dint y),
\]
where $q(t,y,x)$ are the Brownian transition probabilities. Since $q$ is solution of the heat equation, so does $h$:
\begin{equation}\label{eq:edp1}
\frac{\partial}{\partial t}\, h(t,x)=\frac{1}{2}\, \frac{\partial ^2}{\partial x^2}\, h(t,x),\quad \forall (t,x)\in\mathbb{R}_+\times\mathbb{R}.
\end{equation}
If $x=\psi(t)$ denotes the unique solution (see  \cite{lerche_1986}) of the implicit equation $h(t,x)=0$, then, defining $C=\{(t,x):\, x\le\psi(t)\}$ and $u(t,x)\dint x:=\mathbb{P}(\tau_\psi>t,\,B_t\in\dint x) $, both $u$ and $h$ satisfy \eqref{eq:edp1} with the particular boundary condition $h(t,\psi(t))=0$. Using an uniqueness argument, Lerche \cite{lerche_1986} proved that
\[
\mathbb{P}(\tau_\psi>t,\,B_t\in\dint x)=h(t,x)\dint x \quad \mbox{for}\quad (t,x)\in C.
\]
Lerche presented also another proof of this result by using martingales. The distribution of the hitting time can therefore be deduced:
\begin{equation}\label{eq:pdf-lerche}
p_\psi(t)=-\frac{d}{dt}\,\left(\int_{-\infty}^{\psi(t)}h(t,x)\dint x\right).
\end{equation}
Furthermore, we can exhibit a general formula which allows to link the function $b(t)$ introduced in \eqref{eq:relat} and \eqref{eq:defdeb} with the measure $F$.
%\[
%b(t)=\frac{1}{2t}\frac{\int_0^\infty yq(t,y,\psi(t))F(\dint y)}{\int_0^\infty q(t,y,\psi(t))F(\dint y)}.
%\] 
The main challenge is then to find appropriate measures $F$ such that $h(t,x)$, $\psi(t)$ and finally $p_\psi(t)$ (or equivalently $b(t)$) are explicit ! Lerche listed few examples (mainly two-sided curved boundaries) containing obviously the straight line case (the corresponding measure $F$ is a Dirac mass).  For instance, $F(\dint y)=\alpha\delta_{c}(\dint y)+(1-\alpha)\delta_{2c}(\dint y)$, with $0<\alpha<1$ and $c$ a real positive number.\\ %$F(\dint y)=\frac{1}{2}\delta_{2}+\frac{1}{2}\delta_4$. 
%then 
%\[
%\psi(t)=1-\frac{t}{2}\,\left(\ln(1+\sqrt{1+8a e^{-4/t}}) -\ln(4a)\right),
%\]
%and 
%\[
%b(t)=\frac{1}{2t}\frac{q(t,2,\psi(t))+2q(t,4,\psi(t))}{q(t,2,\psi(t))+q(t,4,\psi(t))}.
%\]
In order to conclude the Brownian hitting time study, let us mention that only few situations permit to compute explicitly the pdf and to simulate easily the corresponding stopping time.
%%%%%%%%%%%%%%%%%%%%%%%%%%
\subsection{The Bessel case}
%%%%%%%%%%%%%%%%%%%%%%%%%%
 This section aims at introducing properties concerning the Bessel process hitting times. 
%We have previously introduced results concerning one sided boundaries in the Brownian case. 
The link with the Brownian motion study is the following: the Euclidean norm of a $\delta$-dimensional Brownian motion is a $\delta$-dimensional Bessel process denoted by $(X^{\delta,x}_t,\, t\ge 0)$, where $x$ is the starting point. The time needed by the Brownian motion to exit from a sphere of radius $L$ and the passage time through the level $L$ for the Bessel process are therefore identical in distribution. Thus, it would be of prime interest to obtain explicit expressions of the Bessel hitting time pdf.\\
Let us consider first a constant function $\psi(t)=L>0$. The hitting time is defined by
\[
\tau_\psi=\inf\{ t\ge 0:\ X_t^{\delta,x}=\psi(t) \}.
\]
An explicit form of the Laplace transform is available (see \eqref{laplacebess}) and
%can be proved:
%\[
%\mathbb{E}_{x_0}[e^{-\lambda\tau_\psi}]=\frac{(x_0)^{-\nu}}{L^{-\nu}}\frac{I_\nu(x_0\sqrt{2\lambda})}{I_\nu(L\sqrt{2\lambda})},\quad x_0>0,
%\]
%where $I_\nu$ denotes the modified Bessel function and $\nu=\frac{\delta}{2}-1$ is the so-called index of the Bessel process. Ciesielsky and Taylor \cite{Ciesielski-Taylor} 
can be inverted: the expression of the tail distribution involves Bessel functions of the first kind and their positive zeros. Even if the formula is explicit, it is difficult to handle for numerical procedures. To summarise, a simple expression of the pdf in the constant boundary case  is not available. In addition, there is no hope to obtain interesting results for the straight line case and therefore to use tangent approximation for the simulation of general boundaries hitting times !\\
 The only tool which can be helpful in the Bessel case is the method of images. The idea of this procedure has already been presented in the previous section. Let us denote $q_\delta(t,y,x)$ the transition probabilities associated to the Bessel process of dimension $\delta\in\mathbb{N}$ with $\delta> 1$ and let $F$ be a positive $\sigma$-finite measure on $\mathbb{R}_+$. Then 
 \[
h_\delta(t,x):=q_\delta(t,0,x)-\frac{1}{a}\int_0^\infty q_\delta(t,y,x)F(\dint y)
\]
is solution of the following partial differential equation:
\begin{equation}\label{eq:edp}
\frac{\partial}{\partial t}\, h_\delta(t,x)=\frac{1}{2}\, \frac{\partial ^2}{\partial x^2}\, h_\delta(t,x)-\frac{\delta-1}{2}\frac{\partial}{\partial x}\left(\frac{1}{x}h_\delta(t,x)\right),\quad \forall (t,x)\in\mathbb{R}_+\times\mathbb{R}.
\end{equation}
In particular, if $x=\psi(t)$ is defined as the unique solution of $h_\delta(t,x)=0$ (see \cite{deaconu-herrmann13}), then the density of the measure $\mathbb{P}_0(\tau_\psi>t, \, X^{\delta,0}_t\in\dint x)$ and $h_\delta(t,x)$ satisfy the same PDE with the same boundary conditions. By uniqueness, we deduce that the Bessel hitting time pdf is given by
\[
p_{\delta,\psi}(t)=-\frac{d}{dt}\left( \int_0^{\psi(t)}h_\delta(t,x)\dint x \right),\quad t>0.
\]
%This proof has also been adapted to Bessel processes with non-integer dimension $\delta$ using arguments from the martingale theory \cite{deaconu-herrmann13-2}. 
It suffices now to find suitable measures $F$ such that $p_{\delta,\psi}$ and $\psi$ are explicit. This is namely the case for $F(\dint y)=y^{2\nu+1}1_{\{y>0 \}}\dint y$. In this case
% \marginpar{figure ?}
\begin{equation}
\label{eq:thechoice}
\psi(t)=\sqrt{2t\ln\frac{a}{\Gamma(\nu+1)t^{\nu+1}2^\nu}}\quad\mbox{and}\quad p_{\delta,\psi}(t)=\frac{1}{2at}\,\psi^{2\nu+2}(t),
\end{equation}
where $\nu$ is the index of the Bessel process. The main feature of this discussion is that, $\tau_\psi$ can be numerically sampled in a very easy way ! We have
\[
\tau_\psi\overset{\Delta}{=}\left(\frac{a}{\Gamma(\nu+1)2^\nu}\right)^{\frac{1}{\nu+1}} e^{-Z},
\]
where $Z$ is Gamma distributed with parameters $\nu+2$ and $\frac{1}{\nu+1}$ (see Proposition A.1 in \cite{deaconu-herrmann13}). In particular for $\delta=2$, $e^{-Z}$ is given by the product of two independent standard uniformly distributed random variables $U_1U_2$.

 Other curved boundaries are available but in the sequel we use only this particular example. To sum up, whereas the method of images brought few nice examples of curved boundaries with explicit Brownian hitting time pdf, this method plays in the Bessel case a central role with the example \eqref{eq:thechoice}.
\section{Construction of the algorithm} 
The aim of this section is to combine the random walk on spheres introduced in Section \ref{sec:muller} with the method of images for Bessel processes developed in Section \ref{sec:hitting} in order to construct an efficient algorithm for the simulation of both the exit time and exit position of the Brownian motion. 

The domain that the Brownian motion of dimension $\delta$ has to exit from is a sphere 
 centred in $0$ and of radius $L$ denoted by $\mathcal{D}$. The starting position of the Brownian motion is $B^\delta_0=x_0$.
 
 The structure of the algorithm is the following: we construct a Markov chain $(X(n))_{n\ge 0}$ which represents the Brownian motion position at random times $(T(n))_{n\ge 0}$. 
 
 \medskip

\noindent 1. The initialisation parameters are $X(0)=B_0^\delta=x_0$ and $T(0)=0$. 

\medskip

\noindent 2. The first step evaluates the exit time and the exit position of the Brownian motion for a \emph{moving sphere} centred in $x_0$. The radius of the moving sphere varies continuously on time: it is equal to $\psi(t)$ given by \eqref{eq:thechoice}. The first exit time of the sphere is the first hitting time of the curved boundary $\psi$ for the Bessel process of dimension $\delta$, since the norm of the Brownian motion has the same distribution as a Bessel process of dimension $\delta$. The first hitting time $T(1)$ has the same distribution as $\tau_\psi$ in \eqref{eq:thechoice} and can easily be generated. The position $B_{T_1}^\delta$ of the Brownian motion is uniformly distributed on the sphere of radius $\psi(\tau_\psi)$ centred in $x_0$. 

Let us note that the curved boundary $\psi$ depends on a parameter $a>0$ which can be arbitrarily chosen. We choose a suitable value of $a$ such that the moving sphere always stays in $\mathcal{D}$.

\medskip
%the initial sphere centred at $0$ and of radius $L$. 
\noindent 3. The next step of the algorithm starts with $X(1)=B^\delta_{\tau_\psi}=B^\delta_{T(1)}$. We consider then the Brownian exit problem of a new \emph{moving sphere} centred at $X(1)$ and of radius $\psi$ with the corresponding parameter $a$ chosen in such a way that the moving sphere remains in $\mathcal{D}$. 
%the initial domain (sphere centred at $0$ of radius $L$). 
The exit time denoted by $R(2)$ is given by $\tau_\psi$ and the exit position is uniformly distributed on the sphere of radius $\psi(\tau_\psi)$ and centred at $X(1)$. The global time becomes $T(2)=T(1)+R(2)$ and the Markov chain satisfies $X(2)=B_{T(2)}^\delta$. And so on...

\medskip

\noindent 4. The algorithm stops as soon as $\Vert X(n)\Vert \ge L-\varepsilon$ where $\varepsilon$ is a fixed small parameter.\\

The outcome of the algorithm  is $(X(n),T(n))$ that is an approximation of the couple \emph{exit position and exit time}. 
%\fbox{
%\begin{minipage}{11.9cm}
\begin{framed}\emph{
\noindent {\bf\sc Algorithm (A$\delta$).} \\[3pt]
Fix $0<\gamma<1$ and a small parameter $\varepsilon>0$.\\[5pt]
{\bf\rm  Initialisation:} Set $X(0)=x_0$, $T(0)=0$, $R(0)=0$. \\[3pt]
{\bf\rm The $n$-th step:} Let
\[
a_{n-1}=\Big(\gamma^2(L-\Vert X(n-1)\Vert)^2e/(\nu+1)\Big)^{\nu+1}\frac{\Gamma(\nu+1)}{2}.
\] 
While $\Vert X(n-1)\Vert < L-\varepsilon$, choose $U_n$ an uniformly distributed random vector on $[0,1]^{\lfloor \nu\rfloor +2}$, $G_n$ a standard Gaussian random variable and $V_n$ an uniformly distributed random vector on the unit sphere of dimension $\delta$ centred at $x_0$. $U_n$, $G_n$ and $V_n$ are independent. Define the current time
\begin{align*}
\left\{\begin{array}{l}
R(n)=\left( \frac{a_{n-1}}{\Gamma(\nu+1)2^\nu}\ U_n(1)\ldots U_n(\lfloor \nu\rfloor +2) \right)^{\frac{1}{\nu+1}}\exp\left\{ -\frac{\nu-\lfloor \nu\rfloor}{\nu+1}\ G_n^2\right\},\\[12pt]
T(n)=T(n-1)+R(n),
\end{array}\right.\end{align*}
%(the hitting time of the curved boundary $\psi(t)$ given by \eqref{eq:thechoice} parametrised with $a=a_{n-1}$, here it is just the simulation of a function of a Gamma random variable).
%Secondly
and, the current position
$X(n)=X(n-1)+\psi(R(n))V_n.$\\ \\
{\bf\rm Outcome:} The first time $\Vert X(n-1)\Vert \ge L-\varepsilon$, the algorithm  stops and the outcomes are: $X(n)$ and $T(n)$. }
%\end{minipage}}
\end{framed}
 
\noindent The choice of the parameters $a_n$ ensures at each step the moving sphere to belong to the initial domain $\mathcal{D}$. \\
This algorithm is very simple to use. We describe the $2$ dimensional case ($\nu=0$) as the next section, dealing with numerical results, will focus on  this particular situation. The algorithm writes 
\begin{framed}
\emph{
{\sc Algorithm (A2)}\\
{\bf\rm The $n$-th step}: Let
\(
a_{n-1}=\frac{\gamma^2 e}{2}(L-\Vert X(n-1)\Vert)^2
\) and let $(U_n,V_n,W_n)$ be a vector of three independent random variables uniformly distributed on $[0,1]$.
Set $R(n)=a_{n-1}U_nV_n$,
$T(n)=T(n-1)+R(n)$ and
\[
X(n)=X(n-1)+\psi(R(n))\left(\begin{array}{c}
\cos(2\pi W_n)\\
\sin(2\pi W_n)
\end{array}\right)\quad\mbox{with}\quad \psi(t)=\sqrt{2t\ln(a_{n-1}/t)}.
\] 
\vspace{0.5cm}
{\bf\rm Outcome:} The first time $\Vert X(n-1)\Vert \ge L-\varepsilon$, the algorithm stops and the outcomes are:
 $X(n)$ and $T(n)$.
}
\end{framed}
In dimension two each step {\bf  only requires to sample three uniform random variables !}

Let us denote by $N_\varepsilon$ the number of steps of the algorithm (A$\delta$). We obtain (see \cite{deaconu-herrmann13}) the following convergence results:
\begin{thm}
\label{thm}
1. There exist $C_\delta>0$ and $\varepsilon_0(\delta)>0$ such that
\begin{equation}\label{eq:rate}
\mathbb{E}[N_\varepsilon]\le C_\delta |\ln(\varepsilon)|, \quad \mbox{for any}\ \varepsilon\le \varepsilon_0(\delta).
\end{equation}
2. As $\varepsilon$ goes to zero, the couple $(X(N_\varepsilon),T(N_\varepsilon))$ converges in probability towards the couple $(B_\tau^{\delta,x_0},\tau)$ where $\tau$ is the Brownian exit time of the sphere centred in $0$ of radius $L$.  
\end{thm}
{\sc Sketch of proof of 2.} The proof of the convergence rate %of the algorithm
%\eqref{eq:rate} 
is based on the potential theory for Markov chains. We present here only the main ideas for the second result of this theorem.
For $\eta>0$, let us prove that 
\begin{equation}\label{eq:conv}
\lim_{\varepsilon \to 0}\mathbb{P}(\{\Vert B^{\delta,x_0}_\tau-X(N_\varepsilon)\Vert>\eta\}\cup\{ |\tau-T(N_\varepsilon)|>\varepsilon\eta \})=0.
\end{equation}
It is easy to obtain that $T(N_\varepsilon)$ converges to $\tau$: the algorithm stops with a Brownian position $X(N_\varepsilon)$ in a $\varepsilon$ neighbourhood of the sphere of radius $L$. Moreover the Brownian path does not hit the sphere before (we deduce that $\tau\ge T(N_\varepsilon)$). Thus, due to the strong Markov property and the rotational invariance of the Brownian motion, we can consider that the paths after time $T(N_\varepsilon)$ have the same behaviour as a Brownian motion starting from $0$,  at a distance less than $\varepsilon$  with respect to some convex surface. So the projection $\overline{B_t}$ of the Brownian motion in the direction corresponding to the minimal distance between the origin and the surface, is then a one-dimensional Brownian motion and we get
\[
A_1^\varepsilon:=\mathbb{P}(|\tau- T(N_\varepsilon)|>\varepsilon\eta)\le \mathbb{P}_0\Big(\sup_{0\le t\le \varepsilon\eta}\overline{B}_t<\varepsilon\Big)\le \sqrt{\frac{2\varepsilon}{\eta \pi}}.
\] 
Moreover,
\begin{align*}
A^\varepsilon_2 &:=\mathbb{P}\Big(\{\Vert B^{\delta,x_0}_\tau-X(N_\varepsilon)\Vert >\eta\}\cap\{ |\tau-T(N_\varepsilon)|\le\varepsilon\eta \}\Big)\\
&\le  \mathbb{P}\Big(\{\sup_{T(N_{\varepsilon} )\le t\le T(N_{\varepsilon}) +  \varepsilon\eta}\Vert B_t^{\delta,x_0}-X(N_\varepsilon)\Vert>\eta\} \cap \{ |\tau-T(N_\varepsilon)| \le\varepsilon\eta \}\ \Big)\\
&\le \mathbb{P}\Big(\sup_{T(N_{\varepsilon} )\le t\le T(N_{\varepsilon}) +  \varepsilon\eta}\Vert \ B^{\delta,x_0}_t-X(N_\varepsilon)\Vert>\eta\Big)
\le \mathbb{P}\Big(\sup_{0\le t\le \varepsilon\eta}\Vert B_t^{\delta,0}\Vert>\eta\Big)\\
&\le 2\delta\, \mathbb{P}\Big(\sup_{0\le t\le \varepsilon\eta}\overline{B}_t \ge \frac{\eta}{\sqrt{\delta}} \Big)=2\delta\,\mathbb{P}\Big(|G|\ge \frac{\sqrt{\eta}}{\sqrt{\varepsilon\delta}}\Big)\to 0,\quad \mbox{as}\quad \varepsilon\to 0.
\end{align*}
Here $G$ stands for a standard Gaussian distributed random variable. Combining the convergence of $A_1^\varepsilon$ and $A_2^\varepsilon$, as $\varepsilon$ goes to 0, leads to \eqref{eq:conv} and finally to the second statement of Theorem \eqref{thm}.\hfill{$\Box$}
\vspace{1cm}
%%%%%%%%%%%%%%%%%%%%%%%%%%%%%%%%
\section{Numerical results}
%%%%%%%%%%%%%%%%%%%%%%%%%%%%%%%%
If we are only interested in boundary valued problem like the Laplace
equation, the WOS method is preferable to the WOMS method because
its number of steps before absorption is obviously smaller. However,
in many situations like the approximation of the leading Laplace operator eigenvalue the simulation of the law of the Brownian exit time from a domain is also required. The WOMS provides this exit
time very naturally in any dimension whereas its simulation is more
difficult using the WOS. Nevertheless the simulation of the exit time $\tau_{r}$
of a sphere of radius $r$ is obtained by $r^{2}\tau_{1}$ (with another starting point) thanks
to scaling arguments.

As a consequence, we just have to sample from
$\tau_{1}$ which can be done by at least two methods. The first one
relies on the inversion method applied to $F(t)=P(\tau_{1}<t)$ written
as its spectral expansion as described in section two. To perform
the inversion method, we compute $F^{-1}(U)$ using Newton's method.
Depending on the value of $U,$ the initialisation of Newton's method
and the number of terms kept in the truncation of the series, need
to be adapted for the method to be efficient. The second one uses
a precomputation of sample values of $\tau_{1}$ stored in a large
file. The idea is to pick uniformly at random one value for the time simulation
in the precomputed
file whenever needed. 

The file is built using the inversion method
or a simulation based on the corrected Euler scheme with a very small
step. 
Nevertheless, the precomputation and the inversion method depend
on the dimension of the sphere and are quite consuming tasks. We focus
on a two dimensional example in order to compare the different
approaches. Our test case is the problem of computing the mean exit
time of the unit circle starting at a given point $(x,y).$ Its exact
value is 
$$\mathbb{E}(\tau_{(x,y)})=\frac{1-(x^{2}+y^{2})}{2}.$$
\subsection{Mean time to absorption }
We first study the mean number of steps $\mathbb{E}(N_{\varepsilon}^{WOS})$
and $\mathbb{E}(N_{\varepsilon}^{WOMS})$ before absorption of the WOS and
WOMS methods respectively as a function of the absorption parameter
$\varepsilon$ (with $\gamma = 0.99$). 

We know from our previous results and from known results
on the WOS that both methods behave like $a+b\left|\ln(\varepsilon)\right|$
for $\varepsilon$ small enough. For the starting point $(0.5,0),$
we plot in the next figure both quantities as a function of $\left|\ln(\varepsilon)\right|$
for $\varepsilon=10^{-n},2\leq n\leq 8$ as well as their least-square
fitting (l.l.sq.) which are respectively equal to
\[
E(N_{\varepsilon}^{WOS})\simeq0.3+1.44\left|\ln(\varepsilon)\right|,E(N_{\varepsilon}^{WOMS})\simeq-3.84+3.41\left|\ln(\varepsilon)\right|.
\]
\centerline{\includegraphics[scale=0.7]{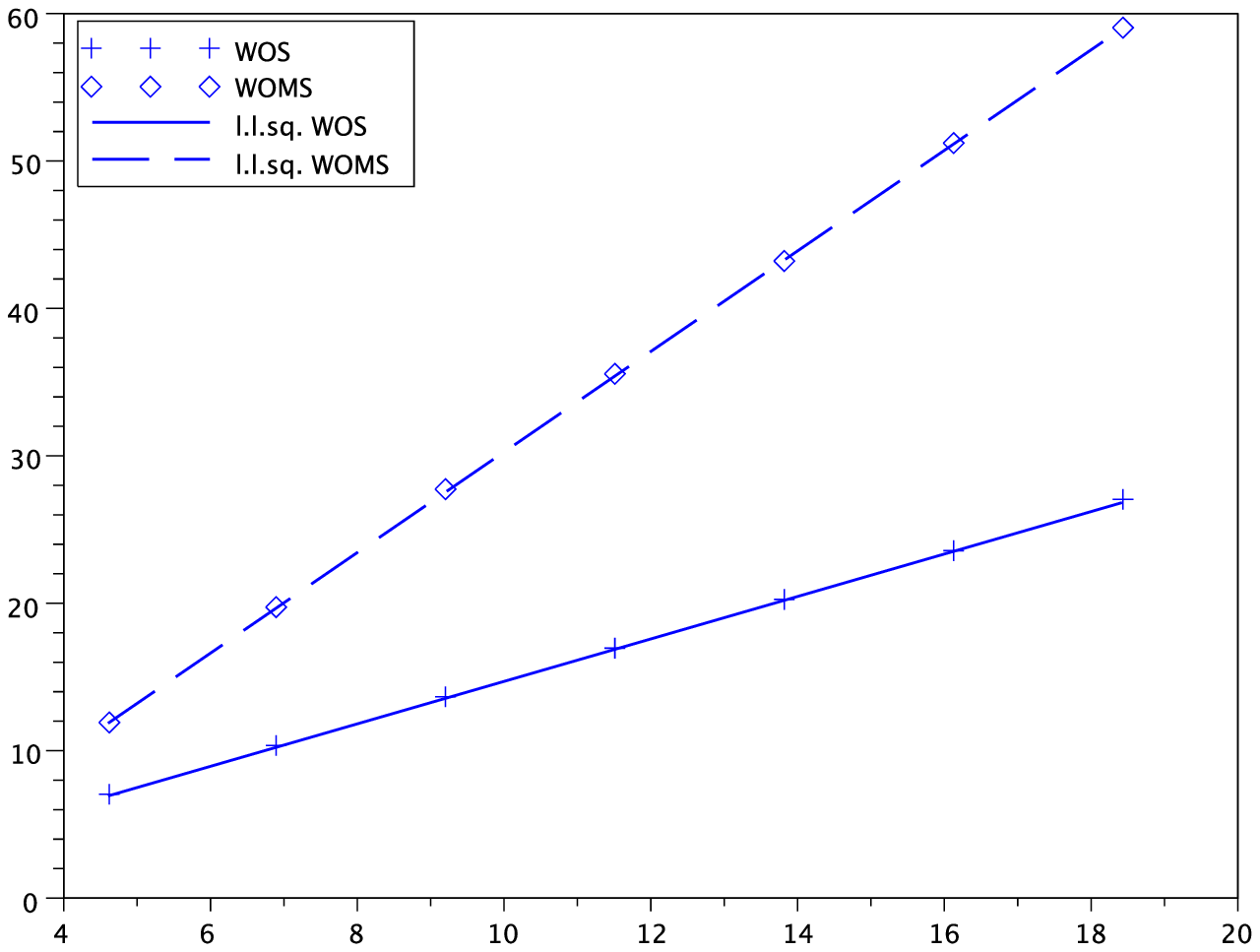}}
\centerline{Fig.1. Average number of steps versus $|\ln(\eps)|$}
\par\null
\vspace*{0.5cm}
\noindent These values have been computed using $10^6$ trajectories. Both
methods fit the model very well and we observe that the number of steps
is around twice bigger for the WOMS method. 

\subsection{Efficiency of the different approaches}

Now we want to study the efficiency in terms of computational times
of the three different approaches for the exit
time of the sphere. The size of the precomputed file is $10^6$
and it should be stored in a binary format to make its opening time
negligible. We have tested all three methods on different starting
points with $\varepsilon=10^{-5}$ and $10^6$ trajectories. They
all gave approximations of the exact mean value correct up to three
or four digits with similar variances. We just need to compare the
computational times in seconds on a standard computer $T^{F},T^{I}$
and $T^{WOMS}$ of the WOS with the precomputed file, the WOS using
the inversion method and of the WOMS respectively. For the starting
point $(0.5,0),$ we plot in the next figure these quantities as a
function of $\left|\ln(\varepsilon)\right|$ for $\varepsilon=10^{-n},2\leq n\leq 8$
as well as their least-square fitting which are respectively equal
to
\[
T^{F}\simeq0.04+0.4\left|\ln(\varepsilon)\right|,T^{I}\simeq0.98+3.63\left|\ln(\varepsilon)\right|,T^{WOMS}\simeq-0.94+1.1\left|\ln(\varepsilon)\right|.
\]
\centerline{\includegraphics[scale=0.7]{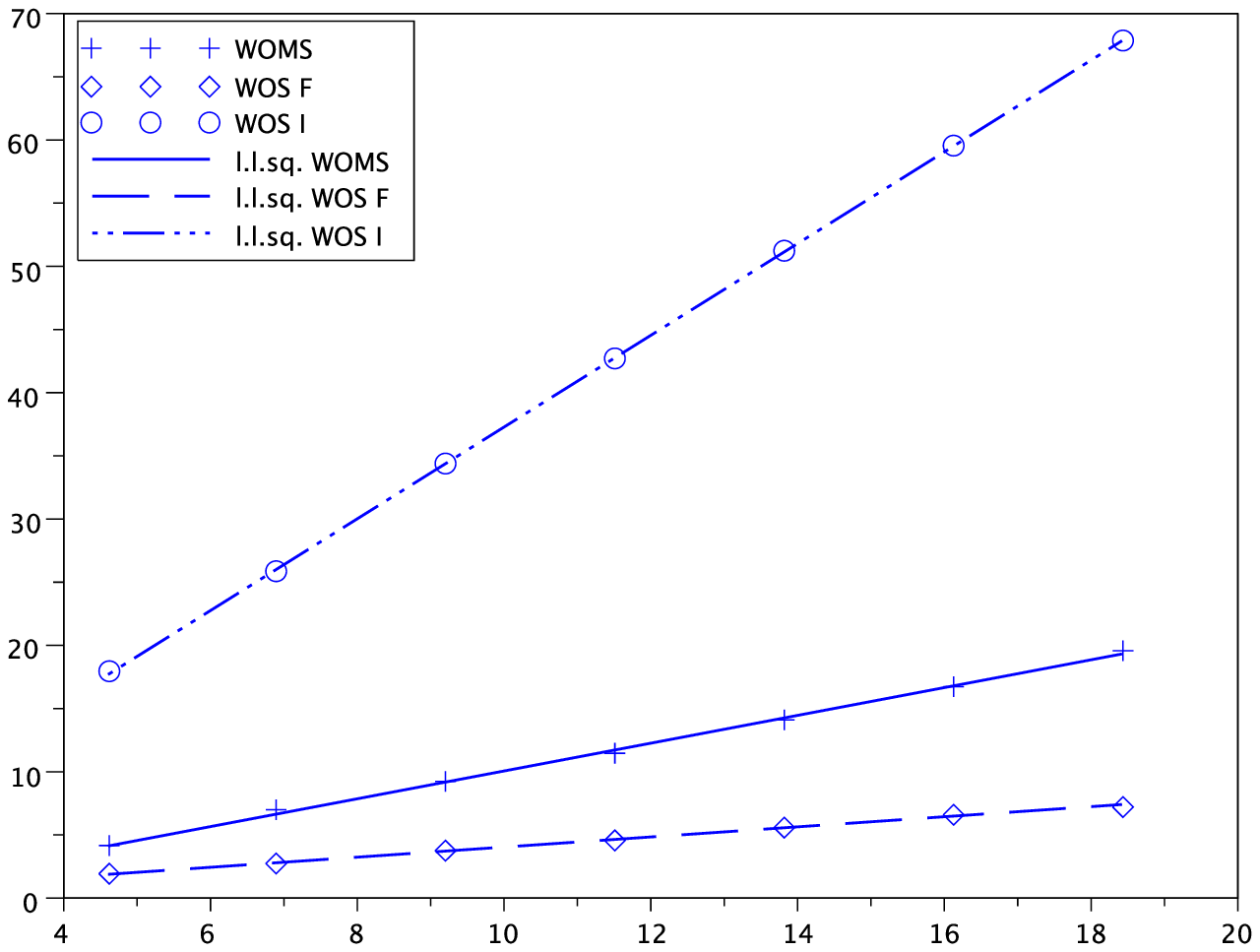}}
\centerline{Fig.2. Cpu versus $|\ln(\eps)|$}

\vspace*{0.2cm}

\noindent Once again, the three methods fit the model very well. We conclude that
the WOMS method is far and away better than the WOS method using
the inversion of the distribution function. It is not surprisingly less
efficient than the WOS coupled with the precomputation. However this
last technique introduces a supplementary bias linked to the size
of the precomputed file which is not easy to quantify and is inextricably linked to a hard precomputed procedure.

%%%%%%%%%%%%%%%%%%%%%%%%%%%%%%%%%
\section{Conclusion}
%%%%%%%%%%%%%%%%%%%%%%%%%%%%%%%%%%

As a conclusion, the walk on moving sphere is a very simple tool to
  compute simultaneously the exit position and the exit time of the
  Brownian motion from a domain in any dimension. It avoids heavy
  computations or additional bias that happened with standard
  techniques. Consequently, we hope the WOMS to replace these
  techniques for applications like principal eigenvalue computations
  \cite{lejay-maire07} where both an accurate and fast simulation of exit times are
  crucial.% which is huge time consuming.

\vspace*{0.2cm}

%\noindent {\bf References}
%\bibliographystyle{elsarticle-num}
%\bibliographystyle{plain}
%\bibliography{htbiblio-1}
\begin{small}

\end{small}
\end{document}